\documentclass[10pt]{article}

\usepackage{amsmath}
\usepackage{amsfonts}
\usepackage{amssymb}

\author{Martin Avenda\~no \thanks{Departamento de Matem\'atica, Facultad de
Ciencias Exactas y Naturales, Universidad de Buenos Aires,  1428
Buenos Aires, Argentina. Research supported by grants CONICET PIP
2461/01 and UBACYT X-112. {\tt mavendar@dm.uba.ar}}}
\date{\today}
\title{The number of real roots \\ of a bivariate polynomial on a line}

\newtheorem{teor}{Theorem}[section]
\newtheorem{coro}[teor]{Corollary}
\newtheorem{lema}[teor]{Lemma}
\newtheorem{obse}[teor]{Remark}
\newtheorem{prop}[teor]{Proposition}
\newtheorem{defi}[teor]{Definition}

\newenvironment{demo}[1][]{{\bf Proof #1: }}{\hfill $\Box$}

\begin{document}

\maketitle

\begin{abstract}
We prove that a polynomial $f\in\mathbb{R}[x,y]$ with $t$ non-zero terms, 
restricted to a real line $y=ax+b$, either has at most $6t-4$ zeroes or 
vanishes over the whole line. As a consequence, we derive an alternative 
algorithm to decide whether a linear polynomial divides a 
sparse polynomial $f\in K[x,y]$ with $t$ terms in 
$[\log(H(f)H(a)H(b))[K:\mathbb{Q}]\log(\deg(f))t]^{O(1)}$ bit operations,
where $K$ is a real number field.
\end{abstract}

\bigskip

\section{Introduction}

The famous Descartes' Rule of Signs, 1641, establishes that the
number of positive real roots of a
polynomial $f\in \mathbb{R}[x]$, counted with multiplicities, is bounded by the number of changes of signs
in its ordered vector of coefficients, disregarding the zeros. As a
direct consequence, the number of different real roots of $f$ is
bounded by $2t-1$, where $t$ is its number of non-zero terms (here
all roots are counted with multiplicities, except $0$ which is counted at 
most once).

\bigskip

There are not yet  natural generalizations of Descartes' Rule of
Signs for the multivariate setting, but a lot of work has been and
is being done for estimating the number of real isolated or
non-degenerate roots (that is where the Jacobian does not vanish,
condition that implies that the root is isolated) of multivariate
square systems of real polynomials in the positive orthant, in terms
of the number of variables and the number of non-zero terms that the
system involves.

\bigskip

The main  result in that direction is due to A.\ Khovanskii
\cite{Kho91}. A simple version of it implies that a square system of $n$
real polynomial equations in $n$ indeterminates, which involves in
total $t$ non-zero terms has at most $(n+1)^t 2^{t(t-1)/2}$
non-degenerate real roots in the positive orthant. Improvements of
Khovanskii's bound have afterwards been obtained by D.\ Perrucci \cite{Per05}
and T.Y.~Li, J.M.~Rojas and X.~Wang \cite{LRW03}, but for general systems
the exponential dependence on the number of non-zero terms $t$ can
not be avoided yet.

\bigskip

In \cite{LRW03}, T.Y.\ Li, M.\ Rojas and X.\ Wang  studied particular
cases of bivariate square systems and showed that the number of
common isolated or non-degenerate roots of a trinomial and a
polynomial with at most $t$ non-zero terms, $t\geqslant 3$,  is bounded by
$2^{t}-2$.

\bigskip

Furthermore, Kushnirenko's Conjecture, formulated in the mid-1970'
(which says that a  square system of $n$ real polynomial equations
in $n$ indeterminates such that the $k$-th polynomial has $t_k$
non-zero terms should have at most $(t_1-1)\cdots (t_n-1)$
non-degenerate roots in the positive orthant) turned out to be
false, by the counter-example provided by B.\ Haas in 2002 for a
system of two trinomials in two variables \cite{Haas02}.

\bigskip

The main result of this article is a refinement of the previous
result for the particular case when the trinomial is a linear
polynomial. Without loss of generality we can assume the linear
polynomial is of the form $y-ax-b$ and we thus study the possible
number of real roots of a bivariate polynomial on a line $y=ax+b$:

\begin{teor}\label{teor_khov}
Let $f=\sum_{i=1}^{t}a_{i}x^{\alpha_{i}}y^{\beta_{i}}\in\mathbb{R}[x,y]$ be
a polynomial with at most $t$ non-zero terms, and  let $a,b\in\mathbb{R}$.
Set $g(x)=f(x,ax+b)$. Then either $g\equiv 0$ or $g$ has at most
$6t-4$ real roots,  counted with multiplicities except for the
possible roots  $0$ and $-b/a$ that are counted at most once.
\end{teor}

To our knowledge, this is the first time a non-exponential bound is
achieved, even for systems of very particular shape like this one.
The tools we use are completely elementary, and we are now studying
the possibility of extending the results for more general systems.

\bigskip

As a consequence of our result we derive an alternative algorithm
for checking if a given linear form $y-ax-b$ divides a polynomial
$f$ in in $K[x,y]$, where $K$ is a real number field.  The number of
bit operations performed by the algorithm  is polynomial in the
degree $[K:\mathbb{Q}]$ of the field extension, in the number $t$ of
non-zero terms of $f$, in the logarithm of the degree of $f$  and in
the logarithmic height of $a$, $b$ and $f$.

\bigskip

The first algorithm for this purpose can be deduced from a more
general result by E.~Kaltofen and P.~Koiran \cite{KaKo05}. They showed
a polynomial-time algorithm for computing all linear factors of a sparse
bivariate polynomial. This result has been further generalized in \cite{AKS06}
and \cite{KaKo06} to an algorithm that computes all the small degree factors
of bi- and multi-variate sparse polynomials.
All these algorithms use a version of the ``gap theorem'' introduced by F.\ 
Cucker, P.\ Koiran and S.\ Smale \cite{CKS99}.
Instead of it, we reduce the problem to the univariate case by considering specializations $f(x,x^{n})$ for small values of $n$.

\section{Proof of Theorem \ref{teor_khov}}

\begin{defi}
Let $f=\sum_{i=0}^{d}a_{i}x^{i}\in\mathbb{R}[x]$ be a non-zero polynomial. We
note by $V(f)$ the number of changes of signs in the ordered vector $(a_d,\dots, a_0)$ of the coefficients of
$f$, disregarding the zeroes. We also set $V(0)=-2$.
\end{defi}

Next result is our crucial ingredient in the proof of Theorem \ref{teor_khov}:

\begin{lema}\label{cambios_de_signo}
Let $f\in\mathbb{R}[x]$. Then $V((x+1)f)\leqslant V(f)$.
\end{lema}
\begin{demo}
We proceed by induction in the number $t$ of non-zero terms of $f$.
The theorem is trivial for $t=0$ and $t=1$.
Now let us suppose that it holds for all $t\leqslant n$.
Let $f\in\mathbb{R}[x]$ with $n+1$ non-zero monomials.
$$f=\sum_{i=1}^{n+1}a_{i}x^{\alpha_{i}}\text{ where } a_{i}\neq 0 
\text{ for all $i$ and } 0\leqslant\alpha_{1}<\alpha_{2}<\cdots<\alpha_{n+1}=d=
\deg(f)$$
Let $g=\sum_{i=1}^{n}a_{i}x^{\alpha_{i}}$.
By inductive hypothesis we have $V((x+1)g)\leqslant V(g)$.
First, we consider the case $\alpha_{n}<d-1$, i.e., when the terms of $(x+1)g$ 
do not overlap with those of $a_{n+1}x^{d}(x+1)$.
There are two posibilities: if $a_{n}a_{n+1}>0$, then $V((x+1)f)=V((x+1)g)\leqslant V(g)=V(f)$, and if $a_{n}a_{n+1}<0$, then $V((x+1)f)=V((x+1)g)+1\leqslant V(g)+1=V(f)$.
In both cases we have $V((x+1)f)\leqslant V(f)$.
Now it only remains the case $\alpha_{n}=d-1$. Here $(x+1)f$ and $(x+1)g$ only differ in their terms of degree $d$ and $d+1$, as shown in the following table.

\begin{center}
\begin{tabular}{|c|c|c|}
\hline
 & $x^{d}$ & $x^{d+1}$ \\
\hline
 $(x+1)g$ & $a_{n}$ & $0$ \\
\hline
 $(x+1)f$ & $a_{n}+a_{n+1}$ & $a_{n+1}$ \\
\hline
\end{tabular}
\end{center}

If $a_{n}a_{n+1}>0$, then $V(f)=V(g)$, and according to the table, we have $V((x+1)f)=V((x+1)g)$.
Therefore $V((x+1)f)\leqslant V(f)$.
On the other hand, if $a_{n}a_{n+1}<0$, then $V(f)=V(g)+1$, but we have three different posibilities for the table, depending whether $|a_{n}|$ is greater, equal or less than $|a_{n+1}|$. Set $s={\rm sgn}(a_{n})$.

\begin{center}
\begin{tabular}{|c|c|c|}
\hline
 & $x^{d}$ & $x^{d+1}$ \\
\hline
 $(x+1)g$ & $s$ & $0$ \\
\hline
 $(x+1)f$ & $s$ & $-s$ \\
\hline
 \multicolumn{3}{c}{Case $|a_{n}|>|a_{n+1}|$} \\
\end{tabular}
\hspace{1cm}
\begin{tabular}{|c|c|c|}
\hline
 & $x^{d}$ & $x^{d+1}$ \\
\hline
 $(x+1)g$ & $s$ & $0$ \\
\hline
 $(x+1)f$ & $0$ & $-s$ \\
\hline
 \multicolumn{3}{c}{Case $|a_{n}|=|a_{n+1}|$} \\
\end{tabular}
\hspace{1cm}
\begin{tabular}{|c|c|c|}
\hline
 & $x^{d}$ & $x^{d+1}$ \\
\hline
 $(x+1)g$ & $s$ & $0$ \\
\hline
 $(x+1)f$ & $-s$ & $-s$ \\
\hline
 \multicolumn{3}{c}{Case $|a_{n}|<|a_{n+1}|$} \\
\end{tabular}
\end{center}

The tables above show that $V((x+1)f)\leqslant V((x+1)g)+1$ for each of the three cases. Using the inductive hypothesis and $V(f)=V(g)+1$, we conclude that $V((x+1)f)\leqslant V(f)$.
\end{demo}

\begin{obse}\label{perturbacion}
Let $f,g\in\mathbb{R}[x]$ and suppose that $g$ has $t$ terms. Then $V(f+g)\leqslant V(f)+2t$.
\end{obse}

Note that the value of $V(0)$ is not relevant for theorem \ref{cambios_de_signo}. The only reason for setting $V(0)=-2$ is the previous remark.

\begin{prop}\label{signos_g}
Let $f\in\mathbb{R}[x,y]$ be a polynomial with $t$ non-zero terms. Then $$V(f(x,x+1))\leqslant 2t-2.$$
\end{prop}
\begin{demo}
We write $f=\sum_{i=1}^{n}a_{i}(x)y^{\alpha_{i}}$, where $0\leqslant\alpha_{1}<\cdots<\alpha_{n}$ and $a_{i}(x)\in\mathbb{R}[x]$, and we set $t_{i}>0$ the number of non-zero terms of $a_{i}$. It is clear that $t=t_{1}+\ldots+t_{n}$.

We define $f_{k}=\sum_{i=k}^{n}a_{i}(x)y^{\alpha_{i}-\alpha_{k}}$ for $k=1,\ldots,n$ and $f_{n+1}=0$.
Lemma \ref{cambios_de_signo} and Remark \ref{perturbacion} imply that the polynomials $f_{k}$ satisfy:
\begin{itemize}
\item $f_{n+1}=0\quad\Longrightarrow\quad V(f_{n+1}(x,x+1))=-2$
\item $f_{k}=y^{\alpha_{k+1}-\alpha_{k}}f_{k+1}+a_{k}(x)\quad\Longrightarrow\quad f_{k}(x,x+1)=(x+1)^{\alpha_{k+1}-\alpha_{k}}f_{k+1}(x,x+1)+a_{k}(x)\quad\Longrightarrow\quad V(f_{k}(x,x+1))\leqslant V(f_{k+1}(x,x+1))+2t_{k}$
\item $f = y^{\alpha_{1}}f_{1}\quad\Longrightarrow\quad f(x,x+1)=(x+1)^{\alpha_{1}}f_{1}(x,x+1)\quad\Longrightarrow\quad V(f(x,x+1))\leqslant V(f_{1}(x,x+1))$.
\end{itemize}
Thus, we conclude that $V(f(x,x+1))\leqslant -2+2(t_{1}+\ldots+t_{n})=2t-2$.
\end{demo}

\bigskip

Before finishing the proof of Theorem \ref{teor_khov} we recall Descartes' Rule of Signs:

\begin{teor}
{\bf (Descartes' rule of signs)} Let $f\in\mathbb{R}[x]$ be a non-zero
polynomial. Then $f$ has at most $V(f)$ positive roots counted with
multiplicities.
\end{teor}

\begin{demo}[of Theorem \ref{teor_khov}]
If $a=0$ or $b=0$, then $g\in\mathbb{R}[x]$ is a polynomial with at most $t$ non-zero terms. Descartes' rule of signs implies that, either $g\equiv 0$ or $g$ has at most $2t-1\leqslant 6t-4$ real roots (counted with multiplicities except for the possible root $0$).
In the case $a\neq 0$ and $b\neq 0$, the real roots of $f(x,ax+b)$ correspond one to one to the roots of $f(bx/a,b(x+1))=\widehat{f}(x,x+1)$, where $\widehat{f}=\sum_{i=1}^{t}a_{i}a^{-\alpha_{i}}b^{\alpha_{i}+\beta_{i}}x^{\alpha_{i}}y^{\beta_{i}}$. Since this bijection preserves the multiplicity of the roots and maps the possible roots $0$ and $-b/a$ of $g$ to the roots $0$ and $-1$ of $\widehat{f}(x,x+1)$, it sufficies to consider the case $a=b=1$, i.e. $g=f(x,x+1)$. Suppose that $g\not\equiv 0$. Descartes' rule of signs and proposition \ref{signos_g} imply that the number of positive roots of $g$ counting with multiplicities is at most $2t-2$. On the other hand, the roots of $g$ in $(-\infty,-1)$ correspond to the positive roots of $0\not\equiv g(-1-x)=f(-1-x,-x)=f_{1}(x,x+1)$, where $f_{1}=\sum_{i=1}^{t}a_{i}(-1)^{\alpha_{i}+\beta_{i}}x^{\beta_{i}}y^{\alpha_{i}}$. Therefore the number of roots (with multiplicities) of $g$ in $(-\infty,-1)$ is also bounded by $2t-2$. Finally, the roots of $g$ in $(-1,0)$ correspond to the positives roots of
$$0\not\equiv (x+1)^{\deg(g)}g\left(\frac{-x}{x+1}\right)=
(x+1)^{\deg(g)}f\left(\frac{-x}{x+1},\frac{1}{x+1}\right)=f_{2}(x,x+1)$$
where $f_{2}=\sum_{i=1}^{t}a_{i}(-1)^{\alpha_{i}}x^{\alpha_{i}}y^{\deg(g)-\alpha_{i}-\beta_{i}}$. Therefore there are at most $2t-2$ of such roots. Taking into account the possible roots $0$ y $-1$, counted each one at most once, we conclude that $g$ has at most $6t-4$ real roots.
\end{demo}

\bigskip

\section{Checking linear factors of a bivariate polynomial}\label{el_algoritmo}

\begin{prop}
Let
$f=\sum_{i=1}^{t}a_{i}x^{\alpha_{i}}y^{\beta_{i}}\in\mathbb{R}[x,y]$.
Let $a,b\in\mathbb{R}$ such that $|b|\neq |1-a|$. Then $y-ax-b\,
|\, f\;\Leftrightarrow\;x^{n}-ax-b\, |\, f(x,x^{n})$ for at least
$6t-3$ odd integers $n\geqslant 3$.
\end{prop}
\begin{demo}
$(\Leftarrow):$ Let $3\leqslant n_{1}<n_{2}<\cdots<n_{6t-3}$ the
$6t-3$ odd numbers for which $x^{n}-ax-b\, |\, f(x,x^{n})$. Let
$w_{i}\in\mathbb{R}$ be a root of $x^{n_{i}}-ax-b$ for each
$1\leqslant i\leqslant 6t-3$. Then
$f(w_{i},aw_{i}+b)=f(w_{i},w_{i}^{n_{i}})=0$ for all $1\leqslant
i\leqslant 6t-3$. This means that $f(x,ax+b)$ has at least $6t-3$ real
roots. Applying theorem~\ref{teor_khov} we conclude that
$f(x,ax+b)\equiv 0$, or simply $y-ax-b\, |\, f$. It only remains
to proof that $w_{i}\neq w_{j}$ for all $i\neq j$. Actually, if
$x^{n_{i}}-ax-b$ and $x^{n_{j}}-ax-b$ had a common root
$w=w_{i}=w_{j}\in\mathbb{R}$, then $w^{n_{i}-n_{j}}=1$ and
therefore $w=\pm 1$. This would imply that $0=w^{n_{i}}-aw-b=-b\pm
(1-a)$, in contradiction with the hypothesis $|b|\neq |1-a|$.
\end{demo}

\begin{coro}
If $f\in\mathbb{R}[x,y]$ has $t$ non-zero terms, then there is and odd integer $3\leqslant n\leqslant 12t-5$ such that $f(x,x^{n})\neq 0$.
\end{coro}
\begin{demo}
Otherwise every polynomial $y-ax-b\in\mathbb{R}[x,y]$ with $|b|\neq|1-a|$ would divide $f$.
\end{demo}

\bigskip

Note that if $(a,b)\neq(0,\pm 1)$, then either $|b|\neq |1-a|$ or
$|b|\neq |1+a|$.

\bigskip

\begin{center}
{\bf Algorithm TEST}
\end{center}

{\bf Input :} A sparse polynomial
$f=\sum_{i=1}^{t}a_{i}x^{\alpha_{i}}y^{\beta_{i}}\in K[x,y]$ with
$t$ monomials, encoded as a list of vectors
$(a_{i},\alpha_{i},\beta_{i})\in
K\times\mathbb{N}_{0}\times\mathbb{N}_{0}$ representing the
monomials of $f$, and two numbers $a,b\in K$.
\smallskip

{\bf Output :} True or False depending whether $y-ax-b\, |\,
f(x,y)$ or not.
\smallskip

{\bf Step 1 :} If $(a,b)=(0,\pm 1)$, compute $f(x,b)$. If this
polynomial is zero, return True. Otherwise return False.
\smallskip

{\bf Step 2 :} If $|b|=|1-a|$ then replace $f$ by $f(-x,y)$ and
$a$ by $-a$.
\smallskip

{\bf Step 3 :} For $n=3,5,7,\ldots,12t-5$ do
\smallskip

{\bf Step 3.1 :} If $f(x,x^{n})\neq 0$ then
\smallskip

{\bf Step 3.1.1 :} Compute all the irreducible factors (with
multiplicities) of $x^{n}-ax-b$ in $K[x]$ using a univariate dense
factorization algorithm.
\smallskip

{\bf Step 3.1.2 :} Compute all the irreducible factors (with
multiplicities) of $f(x,x^{n})$ in $K[x]$ with degree $\leqslant
n$ using a univariate sparse factorization algorithm.
\smallskip

{\bf Step 3.1.3 :} If there is an irreducible factor in the first
list that either does not belong to the second list or belongs but
with less multiplicity, then return False.
\smallskip

{\bf Step 4 :} Return True.

\bigskip

The correctness of the algorithm is a consequence of the previous results. In order to estimate its complexity, we first state the following two famous results on the factorization of polynomials of univariate polynomials.

\begin{description}
\item[DenseFactor] Given $f\in K[x]$ of degree $d$ and absolute height $H$, it is posible to compute all its irreducible factors in $K[x]$ with multiplicities in $[d\, [K:\mathbb{Q}]\, \log H]^{O(1)}$ bit operations (see \cite{LLL82} for the rational case and \cite{SLa85} for the general case).
\item[SparseFactor] Given $f\in K[x]$ a sparse polynomial of degree $d$, with at most $t$ monomials and absolute height $H$, it is posible to find all its irreducible factors (with multiplicities) in $K[x]$ of degree bounded by $s$ in $[t\, s\, [K:\mathbb{Q}]\, \log d\,\log H]^{O(1)}$ bit operations (see \cite{Len99b}).
\end{description}

The complexity of algorithm TEST is clearly dominated by its main loop (step $3$), where it performs $6t-3$ calls to DenseFactor and SparseFactor to factorize $x^{n}-ax-b$ completely and find all the factors of degree bounded by $n$ of $f(x,x^{n})$. We have that $\deg(x^{n}-ax-b)=n\leqslant 12t-5$ and $H(x^{n}-ax-b)\leqslant H(a)H(b)$, therefore step $3.1.1$ requieres at most $[(6t-3)\, (12t-5)\, [K:\mathbb{Q}]\, \log (H(a)\, H(b))]^{O(1)}$ bit operations. On the other hand, we have that $f(x,x^{n})$ is a sparse polynomial with at most $t$ non-zero terms, of degree bounded by $nd\leqslant (12t-5)d$ and absolute height bounded by $(2H)^{t}$ because the coefficients of $f(x,x^{n})$ are sums of at most $t$ coefficients of $f$. Thus, step $3.1.2$ requieres no more than $[(6t-3)\, t\, (12t-5) [K:\mathbb{Q}]\, \log(d(12t-5))\, \log (2H)^{t}]^{O(1)}$ bit operations. This proves that the total number of bit operations performed by the algorithm is polynomial in $t$, $\log(d)$, $[K:\mathbb{Q}]$ and $\log(H(a)H(b)H)$.

\section*{Acknowledgements}
The author thanks Teresa Krick for helping him write this paper, and Daniel
Perrucci for reading an earlier version of this work and for several useful
discussions on fewnomial systems. Thanks also to J.\ Maurice Rojas and 
Frank Sottile for useful discussions.


\begin{thebibliography}{99}
\bibitem{LLL82}
A.K.\ Lenstra, H.W.\ Lenstra, L.\ Lovasz: Factoring polynomials with
rational coefficients. Math.\ Ann.\ 261 (1982), 515-534.
\bibitem{SLa85}
S.\ Landau: Factoring polynomials over algebraic number fields.
SIAM J.\ Comput.\ 14 (1985), 184-195.
\bibitem{Len87} A.K.\ Lenstra: Factoring multivariate polynomials over algebraic number
fields. SIAM J.\ Comput.\ Vol.\ 16 No.\ 3 (1987), 591--598.
\bibitem{Kho91} A.\ Khovanskii: Fewnomials. AMS press, Providence, Rhode Island (1991).
\bibitem{CKS99} F.\ Cucker, P.\ Koiran, S.\ Smale: A polynomial time algorithm for diophantine equations in one
variable. JSC Vol.\ 27 No.\ 1 (1999), 21--29.
\bibitem{Len99b} H.W.\ Lenstra: Finding small degree factors of lacunary polynomials.
Number theory in progress, Vol.\ 1 (1999), 267--276.
\bibitem{Len99a} H.W.\ Lenstra: On the factorization of lacunary polynomials. Number theory
in progress, Vol.\ 1 (1999), 277--291.
\bibitem{Haas02} B.\ Haas: A simple counter-example to Kouchnirenko's conjecture. Beitr\"age zur Algebra und Geometrie, Vol.\ 43 No.\ 1 (2002), 1--8.
\bibitem{LRW03} T.Y.\ Li, J.M.\ Rojas, X.\ Wang: Counting real connected components of trinomial curves intersections and $m$-nomial hypersurfaces. Discrete and computational geometry, Vol.\ 30 No.\ 3 (2003), 379--414.
\bibitem{KaKo05} E.\ Kaltofen, P.\ Koiran: On the complexity of factoring bivariate supersparse
(lacunary) polynomials. ISSAC (2005).
\bibitem{Per05} D.\ Perrucci: Some bounds for the number of connected components of real zero sets of sparse polynomials. Discrete and computational geometry, Vol.\ 34 No.\ 3 (2005), 475--495.
\bibitem{AKS06} M.\ Avenda\~no, T.\ Krick, M.\ Sombra: Factoring
bivariate lacunary polynomials. Journal of Complexity (2006).
\bibitem{KaKo06} E.\ Kaltofen, P.\ Koiran: Finding small degree factors of multivariate supersparse (lacunary) polynomials over algebraic number fields. ISSAC (2006).
\end{thebibliography}
\end{document}